\ifx\documentclass\undefined
   \documentstyle[12pt,a4wide,twoside]{article}
\else
   \documentclass[12pt,twoside]{article}
   \usepackage{a4wide}
\fi

\usepackage{amssymb}
\usepackage{color}
\usepackage{graphicx}
\newcommand	{\nc}	{\newcommand}
\nc	{\be}	{\begin{equation}}
\nc	{\ee}	{\end{equation}}
\nc	{\beq}	{\begin{eqnarray}}
\nc	{\eeq}	{\end{eqnarray}}
\nc	{\beqs}	{\begin{eqnarray*}}
\nc	{\eeqs}	{\end{eqnarray*}}

\nc	{\supp}	{{\rm supp}}
\nc	{\diag}	{{\rm diag}}

\nc{\D}		{\displaystyle}
\nc{\SSS}	{\scriptscriptstyle}

\nc{\comment} [1] {}

\newcount\hour \newcount\minute
\hour=\time    \divide\hour by 60
\minute=\hour  \multiply\minute by 60 \advance\minute by -\time \minute=-\minute
\def\twodigits#1{\ifnum #1 < 10{0#1}\else{#1}\fi}

\renewcommand{\comment}[1]{}

\font\Blackbrd=msbm10 scaled 1200       

\renewcommand{\le}  {\lesssim}      

\renewcommand{\P}   {{\cal P}}      

\newcommand{\R}     {\mbox{\Blackbrd R}}    
\newcommand{\N}     {\mbox{\Blackbrd N}}    
\newcommand{\Poly}  {\mbox{\Blackbrd P}}    

\renewcommand{\div} {{\rm{div} \,}}     

\newtheorem{theorem}{Theorem}[section]
\newtheorem{lemma}[theorem]{Lemma}
\newtheorem{corollary}[theorem]{Corollary}
\newcommand{\bt}{\begin{theorem}}
\newcommand{\et}{\end{theorem}}
\newcommand{\br}{\begin{remark}}
\newcommand{\er}{\end{remark}}
\newcommand{\bc}{\begin{corollary}}
\newcommand{\ec}{\end{corollary}}
\newcommand{\bl}{\begin{lemma}}
\newcommand{\el}{\end{lemma}}
\newcommand{\bp}{\begin{proposition}}
\newcommand{\ep}{\end{proposition}}
\newcommand{\bd}{\begin{definition}}
\newcommand{\ed}{\end{definition}}
\newcommand{\bex}{\begin{example}}
\newcommand{\eex}{\end{example}}

\newtheorem{definition}[theorem]{Definition}
\newtheorem{example}[theorem]{Example}

\newtheorem{remark}[theorem]{Remark}

\newenvironment{proof}{
{\noindent \bf Proof:}}{\quad \hfill \rule{2mm}{2mm}\medskip}

\author{Emmanuel
Creus\'e\footnote{Universit\'e des Sciences et Technologies de
Lille, Laboratoire Paul Painlev\'e UMR 8524, and EPI SIMPAF - INRIA
Lille Nord Europe,
Cit\'e Scientifique, 59655 Villeneuve d'Ascq Cedex email:
creuse@math.univ-lille1.fr}, Serge Nicaise\footnote{Universit\'e de
Valenciennes et du Hainaut Cambr\'esis, LAMAV, FR CNRS 2956,
Institut des Sciences et Techniques de Valenciennes, F-59313 -
Valenciennes Cedex 9 France, email:
Serge.Nicaise@univ-valenciennes.fr}, Emmanuel Verhille \footnote{Universit\'e des Sciences et Technologies de
Lille, Laboratoire Paul Painlev\'e UMR 8524,
Cit\'e Scientifique, 59655 Villeneuve d'Ascq Cedex email:
verhille@math.univ-lille1.fr} }

\begin{document}
\title{Robust equilibrated a posteriori error estimators for the Reissner-Mindlin system}

  {
\maketitle
 \noindent
{
 \begin{abstract}
{We consider a conforming finite element approximation  of  the
Reissner-Mindlin system. We propose a new robust a posteriori error
estimator
 based on $H( \div)$ conforming finite elements
and equilibrated fluxes. It is shown that this estimator gives rise
to an upper bound where the constant is one up to higher order
terms. Lower bounds can also be established with  constants
depending on the shape regularity of the mesh.  The reliability and
efficiency of the proposed estimator are confirmed by some numerical
tests.}
\end{abstract}


 }
\thispagestyle{plain}
\section{Introduction}
\label{introduction} The finite element method is often used for the
numerical approximation of partial differential equations, see,
e.g., \cite{brenner:94,brezzifortin,ciarlet:78}. In many engineering
applications, adaptive techniques based on a posteriori error
estimators have become an indispensable tool to obtain reliable
results. Nowadays there exists a vast amount of literature on
locally defined a posteriori error estimators for problems in
structural mechanics. We refer to the monographs
\cite{AinsworthOden,BS01,neittaanmaki:04,Ver96} for a good overview
on this topic. In general,  upper and lower bounds are established
in order to guarantee the reliability and the efficiency of the
proposed estimator. Most of the existing approaches  involve
constants depending on the shape regularity of the elements; but
these dependencies are often not given.   Only a small number of
approaches gives rise to estimates with explicit constants, see,
e.g.,
\cite{AinsworthOden,braess:06,DN08,ern_nic:07,Repin06,LL83,LW04,neittaanmaki:04,nicaise-wohlmuth:05}.
However in practical applications the knowledge of such constants is
of great importance, especially for adaptivity.

The finite element approximation of the Reissner-Mindlin system
recently became an active subject of research due to its practical
importance and its non trivial challenges to overcome. In
particular, appropriated finite elements have to be used in order to
avoid shear locking. Such elements are in our days well known and
different a priori error estimates are available in the literature.
On the contrary for a posteriori error analysis only a small number
of results exists, we  refer to
\cite{Beiroetco:08,Carstensen:02,CH08,CarstensenSchoberl:06,Repin06,Liberman:01,LovadinaStenberg:05,HuHuang:10}.
Most of these papers enter in the first category mentioned before
and to our knowledge only the paper \cite{Repin06} proposes an
estimator where an upper bound is proved with a constant 1.
 Hence our goal is to give an estimator that is robust with respect
to the thickness parameter $t$, with an explicit constant in the
upper bound, that is also efficient  and that is explicitly
computable. For these purposes
 we use an approach based  on equilibrated fluxes and $H(\div)$--conforming elements.
 Similar ideas can be found, e.g., in
\cite{braess:06,DN08,Repin06,LW04,nicaise-wohlmuth:05}.
 For an overview on
equilibration techniques, we refer to \cite{AinsworthOden,LL83}.

The outline  of the paper is as follows: We recall, in Section 2,
the Reissner-Mindlin system, its numerical approximation and
introduce some useful quantities. Section 3 is devoted to some
preliminary results in order to prove the upper bound. This one
directly follows from these considerations and is given in
full details in section 4. The lower bound developped in section \ref{efficiency} relies on suitable norm
equivalences and by using appropriated $H(\div)$ approximations of
the solutions.  Finally some numerical tests are presented
 in section \ref{sec-num}, that confirm the reliability and the efficiency of our
error estimator.
\section{The boundary value problem and its discretization}
Let $\Omega$ be a bounded open domain of $\R^2$  with a Lipschitz
boundary $\Gamma$ that we suppose to be polygonal. We consider the
following Reissner-Mindlin problem : Given $g \in L^2(\Omega)$
defined as the scaled transverse loading function and $t$ a
fixed positive real number that represents the thickness of the
plate, find $(\omega , \phi) \in  H_0^1(\Omega) \times
H_0^1(\Omega)^2$ such that
\begin{equation}
\label{S2}
a(\phi, \psi) + (\gamma, \nabla v - \psi) = (g,v) \mbox{ for all }
(v,\psi) \in H_0^1(\Omega) \times H_0^1(\Omega)^2,
\end{equation}
where
 \begin{equation}
 \label{S3}
 \gamma =
\lambda \, t^{-2}  \, (\nabla \omega - \phi)  \mbox{ and } a(\phi,
\psi) = \int_{\Omega} \mathcal{C} \varepsilon (\phi)  \varepsilon
(\psi) dx.
\end{equation}
Here, $(\cdot \, , \, \cdot)$ stands for the usual inner product in (any power of) $L^2(\Omega)$,   the operator $:$ denotes the usual term-by term tensor product and
$$\varepsilon(\phi) = \displaystyle \frac{1}{2} (\nabla \phi + {(\nabla \phi)}^T).$$ $\mathcal{C}$ is  the usual elasticity tensor given by
$$
\mathcal{C} \varepsilon (\phi)=2 \, \mu \, \varepsilon
(\phi)+\widetilde{\lambda} \, tr(\varepsilon (\phi)) \, \mathcal{I}.
$$
The parameters $\mu$, $\widetilde{\lambda}$ and $\lambda$  are some
Lam\'e coefficients defined according to the  Young modulus $E$ and
the Poisson coefficient $\nu$ of the material. In the following, for shortness the $L^2(D)$-norm is denoted
by $\|\cdot\|_D$.   The usual norm and seminorm of $H^{1}(D)$ are
respectively denoted by $\|\cdot\|_{1,D}$ and $|\cdot|_{1,D}$ and
the usual norm on $H^{-1}(D)$ is denoted $\|\cdot\|_{-1,D}$.  For
all these norms, in the case $D=\Omega$, the index $\Omega$ is
dropped. The usual Poincar\'e-Friedrichs constant in $\Omega$ is the
smallest positive constant $c_F$ such that
$$
||\phi|| \leq c_F \, |\phi|_1 \quad \forall \phi \in H_0^1(\Omega)^2.
$$
By Korn's inequality \cite{GR86}, $a$ is an inner product on
$H_0^1(\Omega)^2$ equivalent to the usual one. Indeed, defining the
energy norm $||\cdot||_{\mathcal{C}}$ by
$$
\|\psi\|^2_{\mathcal{C}} = a(\psi,\psi) \, \forall \,  \psi \in H_0^1(\Omega)^2,
$$
it can be shown (see annex \ref{annexe1}) that
\begin{equation}
\label{Korn}
|\psi|^2_{1} \leq \displaystyle \frac{1}{\mu} \|\psi\|^2_{\mathcal{C}} \, \forall \,  \psi \in H_0^1(\Omega)^2.
\end{equation}
Consequently, the continuous problem (\ref{S2})-(\ref{S3}) is
well-posed. \bl The problem (\ref{S2})-(\ref{S3}) has a unique
solution $(\omega, \phi) \in H_0^1(\Omega) \times H_0^1(\Omega)^2$.
\el
\begin{proof}
Defining the functional $F((\omega, \phi), (v, \psi)) = a(\phi,\psi)
+ (\gamma, \nabla v - \psi)$ with  $\gamma = \lambda \, t^{-2}  \,
(\nabla \omega - \phi)$, let us  establish its coerciveness, namely
that there exists $k>0$ such that
\begin{equation} \label{coercive}
F((\omega, \phi), [\omega, \phi ]) \geq k \, (|\omega|^2_{1} +
|\phi|^2_{1} ), \forall (\omega, \phi) \in H_0^1(\Omega) \times
H_0^1(\Omega)^2.
\end{equation}
Fix an arbitrary pair $(\omega, \phi) \in H_0^1(\Omega) \times
H_0^1(\Omega)^2.$ First of all, (\ref{Korn}) and the standard
Cauchy-Schwarz inequality lead  to
\begin{eqnarray}
\nonumber F((\omega, \phi), (\omega, \phi)) & \geq& \mu |\phi|^2_{1}
+ \lambda t^{-2} \left((1-\eta) | \omega |^2_{1} +
\left(1-\frac{1}{\eta}\right) \|\phi\|^2 \right), \forall   \eta >
0.
\end{eqnarray}
Then we directly obtain
\begin{eqnarray}
F((\omega, \phi), (\omega, \phi))  \geq  \frac{\mu}{2} |\phi|^2_{1}
+ \lambda t^{-2} (1-\eta) | \omega |^2_{1} \label{rel1} + \left(
\frac{\mu}{2 \, c^2_F} + \lambda t^{-2}
\left(1-\frac{1}{\eta}\right) \right) \|\phi\|^2, \forall   \eta >
0.
\end{eqnarray}
Choosing now $\displaystyle \eta = \frac{2 c^2_F \lambda
t^{-2}}{\mu+2 c^2_F \lambda t^{-2}} <1$ in (\ref{rel1}), we have
$$
F((\omega, \phi), (\omega, \phi)) \geq \frac{\mu}{2} |\phi|^2_{1} +
\frac{\mu \, \lambda t^{-2}}{\mu+2 \, c^2_F \lambda t^{-2}} |\omega
|^2_{1}.
$$
This shows that (\ref{coercive}) holds with $\displaystyle k = \min
\left(\frac{\mu}{2}, \frac{\mu \, \lambda t^{-2}}{\mu +2 \, c^2_F
\lambda t^{-2}}\right)$. The conclusion follows from the Lax-Milgram
lemma for which the other assumptions to fulfill are obvious.
\end{proof} \\

Let us now consider a discretization of (\ref{S2})-(\ref{S3}) based
on a conforming triangulation $\mathcal{T}_h$ of $\Omega$ composed
of triangles. We assume that this triangulation is regular, i.e.,
for any element $T \in \mathcal{T}_h$, the ratio $h_T/\rho_T$ is
bounded by a constant $\sigma>0$ independent of $T$ and of the mesh
size $h= \displaystyle \max_{T \in \mathcal{T}_h} h_T$, where $h_T$
is the diameter of $T$ and $\rho_T$ the diameter of its largest
inscribed ball. We consider on this triangulation the classical
conforming $\Poly_1$ finite element spaces $W_h \times \Theta_h$
defined by
$$
W_h= \left\{ v_h \in {\cal{C}}^0(\bar \Omega) ; v_h=0 \mbox{ on }
\partial \Omega \mbox{ and } {v_h}_{|T} \in \Poly_1 (T) \;  \forall
T \in \mathcal{T}_h \right\} \subset H_0^1(\Omega),
$$
$$
\Theta_h=W_h \times W_h \subset H_0^1(\Omega) \times H_0^1(\Omega).
$$
The discrete formulation of the Reissner-Mindlin problem is now to
find $(\omega _h , \phi_h) \in W_h \times \Theta _h$ such that
 \begin{equation}\label{pbdiscret}
a(\phi_h, \psi_h) + (\gamma_h, \nabla v_h - \textbf{R}_h \psi_h) =
(g,v_h) \mbox{ for all } (v_h,\psi_h) \in W_h \times \Theta_h,
\end{equation}
\noindent with
\begin{equation}\label{gammah}
\gamma_h = \lambda t^{-2} (\nabla \omega_h - \textbf{R}_h \phi_h).
\end{equation}
Here, $\textbf{R}_h$ denotes the reduction integration operator in
the context of shear-locking with values in the so-called discrete
shear force space $\Gamma_h$  which depends on
the finite element involved \cite{ BBF89, BD85, DHHLR03, DL92,
SS97}. 
We assume moreover that
$$
\mbox{For all } \psi_h \in \Theta_h, \textbf{R}_h \psi_h \in H_0(rot,\Omega),
$$
\noindent where $H_0(rot,\Omega)=\{v \in {L^2(\Omega)}^2 ; \, rot\; v
\in L^2(\Omega) \mbox{ and } v \cdot \tau=0 \mbox{ on }
\partial \Omega\}$, equipped with the norm
$$
 \|v\|^2_{H(rot,\Omega)}=\|v\|^2_{\Omega}+\|rot \,  v \|^2_{\Omega}.
$$
Here, for any  $v=(v_1,v_2)^T \in L^2(\Omega)^2$, $\displaystyle \;
rot \,  v=\partial v_2/\partial x-\partial v_1/\partial y$ and
$\tau$ is the unit tangent vector along $\partial \Omega$.
In this work, $\textbf{R}_h$ is defined as the interpolation operator
from $\Theta_h$ on the $H_0(rot,\Omega)$ conforming 
lower-order Nedelec finite element space \cite{GR86}. \\
\\
\noindent
By the usual Helmholtz decomposition of any $H_0(rot,\Omega)$ vector field
\cite[p. 299]{brezzifortin}, there exists $w \in H_0^1(\Omega)$ and
$\beta \in H_0^1(\Omega)^2$ such that~:
\begin{equation}\label{rh}
(\textbf{R}_h-I)\phi_h = \nabla w-\beta,
\end{equation}
\noindent as well as a constant $C>0$ such that
$$
\|w\|_{1}+\|\beta\|_{1} \leq C \; \|(\textbf{R}_h-I)\phi_h\|_{H(rot,\Omega)}.
$$
More precisely, we introduce the constant $c_R$ such that
$$
|\beta|_1 \leq c_R \,
\|rot(\textbf{R}_h-I)\phi_h\|,
$$
which can be evaluated by \cite{GR86}
$$
c_{R} = \left( \inf_{q \in L^2(\Omega)} \sup_{v \in
H_0^1(\Omega)^2}\frac{(div \, v, q)}{\|q\| \;  |v|_1} \right)^{-1}.
$$
 Given the exact solution $(\omega, \phi) \in H_0^1(\Omega) \times
H_0^1(\Omega)^2$ as well as the approximated one $(\omega_h,\phi_h)
\in W_h \times \Theta_h$, the usual error $e_h^{rot}$ is defined as
\begin{equation}
\label{ehrot2}
\displaystyle (e^{rot}_h)^2=|\omega - \omega_h|^2_{1} +
|\phi-\phi_h|^2_{1} + \lambda^{-1} t^2
\|\gamma-\gamma_h\|^2 +\lambda^{-2} t^4
\|rot(\gamma-\gamma_h)\|^2+\|\gamma-\gamma_h\|^2_{-1}.
\end{equation}
The residuals are also defined as follows
\begin{eqnarray}
Res_1 (v) &=& (g,v) - (\gamma _h,\nabla
v) \quad \mbox{for all }v \in H^1_0(\Omega), \label{res1}\\
Res_2 (\psi) &=& -a(\phi _h,\psi) + (\gamma _h,\psi) \quad
\mbox{for all }\psi \in {H^1_0(\Omega)}^2. \label{res2}
\end{eqnarray}
Finally, let us now introduce, in the spirit of \cite{Repin06}, the spaces $N_{div}(\Omega)$ and $H_{div}(\Omega)$ respectively defined by
$$
\left.
\begin{array}{l}
H_{div}(\Omega)=\{y \in L^2(\Omega,\R^2)|\,div
\, y \in L^2(\Omega)\}, \\[6pt]
N_{div}(\Omega)=\{x \in L^2(\Omega,\mathcal{M}^2_S)|\,div
\, x \in L^2(\Omega,\R^2)\},
\end{array}
\right.
$$
where $\mathcal{M}^2_S$ is the space of symmetric tensors of second
rank. We now fix an arbitrary $y^* \in H_{div}(\Omega)$  such that
$div \, y^* =-\Pi_h g$, where $\Pi_h$ is the projection operator
from $L^2(\Omega)$ to the piecewise constant fonctions on the
triangulation. Let us also fix  $x^* \in N_{div}(\Omega)$ such that
$div \, x^* =-\gamma_h$. Their existence and construction will be
explained later on.

We finally need to introduce   the following mesh-dependent norm.
For all $(\psi,v) \in H_0^1(\Omega) \times H_0^1(\Omega)^2$, we
define
\begin{equation}\label{norme1h}
|\|(\psi,v)|\|^2_{1,h}=\|\nabla \psi\|^2 +\sum_{T \in
\mathcal{T}_h}\frac{1}{t^2+h_T^2}\|\nabla v-\psi\|^2_{T}.
\end{equation}
For all functional $F$ defined on $H_0^1(\Omega) \times
H_0^1(\Omega)^2$, the dual norm associated with (\ref{norme1h}) is
classically defined by
\begin{equation}
|\|F|\|_{-1,h}=\sup_{(\psi,v) \in H_0^1(\Omega) \times H_0^1(\Omega)^2 \setminus
\{0\}}\frac{F(\psi,v)}{|\|(\psi,v)|\|_{1,h}}.
\end{equation}
\section{Preliminary results}
The aim of this section is to prove four lemmas which will be used in
the following of the paper.
\bl \label{lemma1} Let us consider $(\alpha, \varepsilon) \in
(\R_+^*)^2$. Then we have
$$
\left.
\begin{array}{l}
\displaystyle \lambda (t^{-2}-\alpha^2) \|\nabla
(\omega-\omega_h)-(\phi-\textbf{R}_h
\phi_h)\|^2+\lambda \alpha^2(1-2\varepsilon)\|\nabla
(\omega-\omega_h)\|^2\\[6pt]
\displaystyle \leq \lambda^{-1} t^2
\|\gamma-\gamma_h\|^2-\lambda \alpha^2
\left(1-\frac{2}{\varepsilon}\right) \|\phi_h-\textbf{R}_h
\phi_h\|^2-\lambda
\alpha^2\left(1-\frac{1}{\varepsilon}-\varepsilon\right)\|\phi-\phi_h\|^2.
\end{array}
\right.
$$
\el
\begin{proof}
We first write
$$
\left.
\begin{array}{l}
\displaystyle \|\nabla
(\omega-\omega_h)-(\phi-\phi_h)-(\phi_h-\textbf{R}_h
\phi_h)\|^2\\[6pt]
\displaystyle=\|\nabla
(\omega-\omega_h)\|^2+\|\phi-\phi_h\|^2+\|\phi_h-\textbf{R}_h \phi_h)\|^2 \\[6pt]
\displaystyle -2(\nabla
(\omega-\omega_h),\phi-\phi_h) -2(\nabla
(\omega-\omega_h),\phi_h-\textbf{R}_h
\phi_h)+2(\phi-\phi_h,\phi_h-\textbf{R}_h
\phi_h).
\end{array}
\right.
$$
Consequently, we have
$$
\left.
\begin{array}{lcl}
\displaystyle \lambda^{-1} t^2
\|\gamma-\gamma_h\|^2&=&
\displaystyle \lambda (t^{-2}-\alpha^2) \|\nabla
(\omega-\omega_h)-(\phi-\textbf{R}_h \phi_h)\|^2\\[6pt]
\displaystyle &+& \lambda \alpha^2(\|\nabla
(\omega-\omega_h)\|^2+\|\phi-\phi_h\|^2+\|\phi_h-\textbf{R}_h
\phi_h\|^2)\\[6pt]
\displaystyle &+&2 \, \lambda
\alpha^2(\phi-\phi_h,\phi_h-\textbf{R}_h
\phi_h)-2\lambda \alpha^2(\nabla
(\omega-\omega_h),\phi-\phi_h)\\[6pt]
\displaystyle &-&2 \, \lambda \alpha^2(\nabla
(\omega-\omega_h),\phi_h-\textbf{R}_h \phi_h).
\end{array}
\right.
$$
Using the three following Young inequalities
$$
\left\{
\begin{array}{lcl}
\displaystyle -2(\phi-\phi_h,\phi_h-\textbf{R}_h
\phi_h) &\leq& \varepsilon
\|\phi-\phi_h\|^2+\displaystyle \frac{1}{\varepsilon}\|\phi_h-\textbf{R}_h
\phi_h\|^2,\\[8pt]
\displaystyle 2(\nabla
(\omega-\omega_h),\phi-\phi_h)&\leq& \varepsilon
\|\nabla
(\omega-\omega_h)\|^2+\displaystyle \frac{1}{\varepsilon}\|\phi-
\phi_h\|^2,\\[8pt]
\displaystyle 2(\nabla (\omega-\omega_h),\phi_h-\textbf{R}_h
\phi_h)&\leq& \varepsilon \|\nabla
(\omega-\omega_h)\|^2+ \displaystyle \frac{1}{\varepsilon}\|\phi_h-\textbf{R}_h
\phi_h\|^2,
\end{array}
\right.
$$
we get
$$
\left.
\begin{array}{l}
\lambda (t^{-2}-\alpha^2) \|\nabla
(\omega-\omega_h)-(\phi-\textbf{R}_h \phi_h)\|^2\\[6pt]
\displaystyle \leq \lambda^{-1} t^2
\|\gamma-\gamma_h\|^2-\lambda \alpha^2\left(\|\nabla
(\omega-\omega_h)\|^2+\|\phi-\phi_h\|^2 + \|\phi_h-\textbf{R}_h
\phi_h\|^2 \right)\\[6pt]
\hspace{5mm} \displaystyle + \lambda \alpha^2 \left( \varepsilon
\|\phi-\phi_h\|^2+\frac{1}{\varepsilon}\|\phi_h-\textbf{R}_h
\phi_h\|^2+\varepsilon \|\nabla
(\omega-\omega_h)\|^2+\frac{1}{\varepsilon}\|\phi-\phi_h
\|^2 \right.\\[6pt]
\hspace{5mm} \left.\displaystyle +\varepsilon \|\nabla
(\omega-\omega_h)\|^2+\frac{1}{\varepsilon}\|\phi_h-\textbf{R}_h
\phi_h\|^2 \right)\\[6pt]
\displaystyle  = \lambda^{-1} t^2
\|\gamma-\gamma_h\|^2-\lambda \alpha^2
\left(1-\frac{2}{\varepsilon}\right) \|\phi_h-\textbf{R}_h
\phi_h\|^2-\lambda
\alpha^2\left(1-\frac{1}{\varepsilon}-\varepsilon \right)\|\phi-\phi_h||^2\\[6pt]
\hspace{5mm} \displaystyle -\lambda \alpha^2(1-2\varepsilon)\|\nabla
(\omega-\omega_h)\|^2.
\end{array}
\right.
$$
This proves the lemma.
\end{proof} \\
\bl \label{lemmegamma} we have
\begin{equation}\label{h-1}
\left.
\begin{array}{ll}
\displaystyle\|\gamma-\gamma_h\|^2_{-1}& \leq
\displaystyle 4 \, (\mu+\tilde{\lambda}) \, \|\phi-\phi_h\|^2_{\mathcal{C}}+2\|Res_2\|^2_{-1}.
\end{array}
\right.
\end{equation}
\el
\begin{proof}
First, it can be shown that for any $\psi \in (H_0^1(\Omega))^2$,
$$
\|\psi\|^{2}_{\mathcal{C}} \leq 2 \, (\mu+\tilde{\lambda}) |\psi|^2_1,
$$
so that
$$
\begin{array}{lcl}
\displaystyle(\gamma-\gamma_h,\psi)&=&\displaystyle a(\phi-\phi_h,\psi)+a(\phi_h,\psi)-(\gamma_h,\psi)\\[8pt]
&=&\displaystyle a(\phi-\phi_h,\psi)-Res_2(\psi)\\[8pt]
&\leq&\displaystyle\|\phi-\phi_h\|_{\mathcal{C}}\|\psi\|_{\mathcal{C}}+\|Res_2\|_{-1}|\psi|_{1}\\[8pt]
&\leq&\displaystyle
\left((2 \, (\mu+\tilde{\lambda}))^{1/2} \|\phi-\phi_h\|_{\mathcal{C}}+\|Res_2\|_{-1}\right)|\psi|_{1}.
\end{array}
$$
\noindent Hence we get
$$
\left.
\begin{array}{ll}
\displaystyle\|\gamma-\gamma_h\|^2_{-1}& \leq
\displaystyle
\left((2 \, (\mu+\tilde{\lambda}))^{1/2} \|\phi-\phi_h\|_{\mathcal{C}}+\|Res_2\|_{-1}\right)^2\\[8pt]
& \leq \displaystyle 4 \, (\mu+\tilde{\lambda}) \,
\|\phi-\phi_h\|^2_{\mathcal{C}}+2 \, \|Res_2\|^2_{-1}.
\end{array}
\right.
$$
\end{proof}
\bl
\label{residu}
$$
\displaystyle \|\phi-\phi_h\|^2_{\mathcal{C}}+\lambda^{-1} t^2
\|\gamma-\gamma_h\|^2 =
 Res_1(\omega-\omega_h+w)+Res_2(\phi-\phi_h+\beta)-a(\phi-\phi_h,\beta),
$$
where $w$ and $\beta$ are given by the Helmholtz decomposition (\ref{rh}).
\el
\begin{proof}
This result is similar to the one given in \cite{CH08}. First,
(\ref{S2}) and (\ref{rh}) lead to
$$
\left.
\begin{array}{ll}
\displaystyle (\gamma-\gamma_h,(\textbf{R}_h-I)
\phi_h)&=(\gamma-\gamma_h, \nabla
w-\beta)\\[6pt]
&\displaystyle=(\gamma, \nabla w)
-(\gamma, \beta)-(\gamma_h, \nabla w-\beta)\\[6pt]
\displaystyle &=(g,w)-a(\phi,\beta)-(\gamma_h,\nabla
w-\beta)\\[6pt]
\displaystyle
&=-a(\phi-\phi_h,\beta)+(g,w)-a(\phi_h,\beta)-(\gamma_h,\nabla
w-\beta).
\end{array}
\right.
$$
A simple calculation shows that
$$
\left.
\begin{array}{l}
\displaystyle \|\phi-\phi_h\|^2_{\mathcal{C}}+\lambda^{-1} t^2
\|\gamma-\gamma_h\|^2\\[6pt]
\displaystyle=a(\phi-\phi_h,\phi-\phi_h)+(\gamma-\gamma_h,(\nabla
\omega-\nabla
\omega_h)-(\phi-\phi_h))+(\gamma-\gamma_h,(\textbf{R}_h-I)
\phi_h)\\[6pt]
\displaystyle =(g,\omega-
\omega_h)-a(\phi_h,\phi-\phi_h)-(\gamma_h,\nabla
(\omega- \omega_h))\\[6pt]
\displaystyle \hspace{6mm} +(\gamma_h,\phi-\phi_h)
-a(\phi-\phi_h,\beta)+(g,w)-a(\phi_h,\beta)-(\gamma_h,\nabla
w-\beta)\\[6pt]
\displaystyle = Res_2(\phi-\phi_h+\beta)+(g,\omega-
\omega_h+w)-(\gamma_h,\nabla (\omega-
\omega_h+w))-a(\phi-\phi_h,\beta)\\[6pt]
\displaystyle =Res_2(\phi-\phi_h+\beta)+Res_1(\omega-
\omega_h+w)-a(\phi-\phi_h,\beta).
\end{array}
\right.
$$
So we get
$$
\|\phi-\phi_h\|^2_{\mathcal{C}}+\lambda^{-1} t^2
\|\gamma-\gamma_h\|^2=Res_1(\omega-\omega_h+w)+Res_2(\phi-\phi_h+\beta)-a(\phi-\phi_h,\beta).
$$
\end{proof}
\bl
\label{lemmetiti}
$$
\left.
\begin{array}{l}
\displaystyle
\frac{1}{2}\|\phi-\phi_h+\beta\|^2_{\mathcal{C}}+\frac{1}{2}\|\phi-\phi_h\|^2_{\mathcal{C}}+\frac{1}{2}\lambda^{-1}
t^2
\|\gamma-\gamma_h\|^2\\[6pt]
\displaystyle \hspace{6mm}+\frac{1}{2}\sum_{T \in \mathcal{T}_h}
\frac{\lambda}{t^2+h_T^2}\|\nabla
(\omega-\omega_h+w)-(\phi-\phi_h+\beta)\|^2_{T}\\[8pt]
\leq \displaystyle
Res_1(\omega-\omega_h+w)+Res_2(\phi-\phi_h+\beta)+\frac{1}{2}\|\beta\|^2_{\mathcal{C}}.
\end{array}
\right.
$$
\el
\begin{proof}
The proof is once again similar to the one in \cite{CH08}. Because
of (\ref{rh}), we first remark that
$$
\gamma-\gamma_h=\lambda t^{-2}(\nabla \omega-\nabla \omega_h-\phi+\phi_h+\nabla w-\beta),
$$
so that we have for all $T \in \mathcal{T}_h$
$$
\|\nabla (\omega-\omega_h+w)-(\phi-\phi_h+\beta)\|^2_{T} \leq
\lambda^{-2} t^4\|\gamma-\gamma_h\|^2_{T}.
$$
Then,
$$
\left.
\begin{array}{l}
\displaystyle
\frac{1}{2}\|\phi-\phi_h+\beta\|^2_{\mathcal{C}}+\frac{1}{2}\|\phi-\phi_h\|^2_{\mathcal{C}}+\frac{1}{2}\lambda^{-1}
t^2
\|\gamma-\gamma_h\|^2\\[6pt]
\displaystyle \hspace{6mm}+\frac{1}{2}\sum_{T \in \mathcal{T}_h}
\frac{\lambda}{t^2+h_T^2}\|\nabla
(\omega-\omega_h+w)-(\phi-\phi_h+\beta)\|^2_{T}\\[8pt]
\displaystyle \leq
\frac{1}{2}\|\phi-\phi_h+\beta\|^2_{\mathcal{C}}+\frac{1}{2}\|\phi-\phi_h\|^2_{\mathcal{C}}+\frac{1}{2}\lambda^{-1}
t^2 \|\gamma-\gamma_h\|^2+\frac{1}{2}\lambda^{-1} t^2
\sum_{T \in \mathcal{T}_h} \|\gamma-\gamma_h\|^2_{T}\\[8pt]
\displaystyle \leq \lambda^{-1} t^2
\|\gamma-\gamma_h\|^2+\frac{1}{2}
a(\phi-\phi_h+\beta,\phi-\phi_h+\beta)+\frac{1}{2}a(\phi-\phi_h,\phi-\phi_h)\\[8pt]
\displaystyle = \lambda^{-1} t^2
\|\gamma-\gamma_h\|^2+\frac{1}{2}\left(\|\phi-\phi_h\|^2_{\mathcal{C}}+2a(\phi-\phi_h,\beta)+\|\beta\|^2_{\mathcal{C}}\right)+\frac{1}{2}\|\phi-\phi_h\|^2_{\mathcal{C}}\\[8pt]
 \displaystyle =
\|\phi-\phi_h\|^2_{\mathcal{C}}+\lambda^{-1} t^2
\|\gamma-\gamma_h\|^2+\frac{1}{2}\|\beta\|^2_{\mathcal{C}}+a(\phi-\phi_h,\beta).
\end{array}
\right.
$$
From lemma \ref{residu}, we get
$$
\left.
\begin{array}{l}
\displaystyle
\frac{1}{2}\|\phi-\phi_h+\beta\|^2_{\mathcal{C}}+\frac{1}{2}\|\phi-\phi_h\|^2_{\mathcal{C}}+\frac{1}{2}\lambda^{-1}
t^2
\|\gamma-\gamma_h\|^2\\[6pt]
\displaystyle \hspace{6mm}+\frac{1}{2}\sum_{T \in \mathcal{T}_h}
\frac{\lambda}{t^2+h_T^2}\|\nabla
(\omega-\omega_h+w)-(\phi-\phi_h+\beta)\|^2_{T}\\[8pt]
\displaystyle  \leq Res_1(\omega-\omega_h+w)+Res_2(\phi-\phi_h+\beta)-a(\phi-\phi_h,\beta)+\frac{1}{2}\|\beta\|^2_{\mathcal{C}}+a(\phi-\phi_h,\beta) \\
\displaystyle
=Res_1(\omega-\omega_h+w)+Res_2(\phi-\phi_h+\beta)+\frac{1}{2}\|\beta\|^2_{\mathcal{C}}.
\end{array}
\right.
$$
\end{proof}
\section{Reliability of the estimator}
\bt \label{theototo} Let us consider $0 < \varepsilon < 1/2$, as
well as two parameters $\nu_1>0$ and $\nu_2>0$. Moreover, let us
define
$$
A(\varepsilon)=\max\left(\displaystyle \frac{3}{\mu}+ c_F^2 \frac{\frac{1}{\varepsilon}+\varepsilon-1}{\mu(1-2\varepsilon)}+4 (\mu + \tilde{\lambda});1+\frac{t^2}{\lambda(1-2\varepsilon)}\right).
$$
Then,
\begin{equation}\label{relavecA}
\left.
\begin{array}{ll}
(e_h^{rot})^2&\leq A_1 |\|Res_1|\|^2_{-1,h}+A_2 \|Res_2\|^2_{-1}+A_3 \|\phi-\phi_h+\beta\|^2_{\mathcal{C}}\\[10pt]
&\displaystyle \hspace{5mm} +A_4 \|\phi_h-\textbf{R}_h\phi_h\|^2_{H(rot,\Omega)}-
\sum_{T \in \mathcal{T}_h} A_5^T\|\nabla
(\omega-\omega_h+w)-(\phi-\phi_h+\beta)\|^2_{T},
\end{array}
\right.
\end{equation}
\noindent with
$$
\left|
\begin{array}{l}
\displaystyle A_1=\nu_1 A(\varepsilon)^2;\\[10pt]
\displaystyle A_2=\nu_2 A(\varepsilon)^2+2;\\[10pt]
\displaystyle A_3=\displaystyle \frac{1}{\mu}\left(\frac{1}{\nu_1}+\frac{1}{\nu_2}\right)-A(\varepsilon);\\[10pt]
\displaystyle A_4=\max \left( \frac{\frac{2}{\varepsilon}-1}{1-2\varepsilon} \; ; \; 2+2A(\varepsilon)(\mu + \tilde{\lambda}) c_R^2 \right);\\[10pt]
\displaystyle A_5^T=\frac{\lambda
A(\varepsilon)}{t^2+h_T^2}-\frac{1}{\nu_1 (t^2+h_T^2)}, \; \forall \; T \in
\mathcal{T}_h.
\end{array}
\right.
$$
\et
\begin{proof}
First of all, by using lemma \ref{lemma1} and the fact that $0 <
\varepsilon < 1/2$, we get
$$
\left.
\begin{array}{ll}
(e_h^{rot})^2 &\leq \displaystyle
\left(\frac{1}{\mu}+ c_F^2 \; \frac{\frac{1}{\varepsilon}+\varepsilon-1}{\mu (1-2\varepsilon)}\right)\|\phi-\phi_h\|^2_{\mathcal{C}}+\left(1+\frac{t^2}{\lambda
(1-2\varepsilon)}\right)\lambda^{-1} t^2
\|\gamma-\gamma_h\|^2\\[10pt]
&\displaystyle
\hspace{5mm}+\frac{\frac{2}{\varepsilon}-1}{1-2\varepsilon}\|\phi_h-\textbf{R}_h\phi_h\|^2+\lambda^{-2} t^4 \|rot(\gamma-\gamma_h)\|^2+\|\gamma-\gamma_h\|^2_{-1}.
\end{array}
\right.
$$
Then, because of lemma \ref{lemmegamma} as well as
 $$
\left.
\begin{array}{l}
 \lambda^{-2}t^4\|rot(\gamma-\gamma_h)\|^2 \leq \displaystyle \frac{2}{\mu}  \, \|\phi-\phi_h\|^2_{\mathcal{C}}+2 \,
 \| rot(\phi_h-\textbf{R}_h \phi_h ) \|^2,
\end{array}
\right.
$$
we obtain
$$
\left.
\begin{array}{ll}
(e_h^{rot})^2 &\leq \displaystyle
\left(\frac{3}{\mu}+c_F^2 \frac{\frac{1}{\varepsilon}+\varepsilon-1}{\mu(1-2\varepsilon)}+4 (\mu + \tilde{\lambda}) \right)\|\phi-\phi_h\|^2_{\mathcal{C}} + 2\|rot(\phi_h-\textbf{R}_h \phi_h)\|^2 \\[10pt]
&\displaystyle \hspace{5mm}+\left(1+\frac{t^2}{\lambda
(1-2\varepsilon)}\right)\lambda^{-1} t^2
\|\gamma-\gamma_h\|^2+\left(\frac{\frac{2}{\varepsilon}-1}{1-2\varepsilon}\right)\|\phi_h-\textbf{R}_h\phi_h\|^2\\[10pt]
&\displaystyle \hspace{5mm}+2\|Res_2\|^2_{-1}.
\end{array}
\right.
$$
By the definition of $A(\varepsilon)$ as well as lemma
\ref{lemmetiti}, we get
$$
\left.
\begin{array}{ll}
(e_h^{rot})^2 &\leq \displaystyle A(\varepsilon) \Big( 2Res_1(\omega-\omega_h+w)+2Res_2(\phi-\phi_h+\beta)+\|\beta\|^2_{\mathcal{C}} \\[10pt]
& \displaystyle \hspace{5mm} -\|\phi-\phi_h+\beta\|^2_{\mathcal{C}}-\sum_{T \in \mathcal{T}_h}
\frac{\lambda}{t^2+h_T^2}\|\nabla
(\omega-\omega_h+w)-(\phi-\phi_h+\beta)\|^2_{T}\Big)\\[10pt]
&\displaystyle
\hspace{5mm}+\left(\frac{\frac{2}{\varepsilon}-1}{1-2\varepsilon}\right)\|\phi_h-\textbf{R}_h\phi_h\|^2 +2\|Res_2\|^2_{-1}+2\|rot(\phi_h-\textbf{R}_h \phi_h)\|^2.
\end{array}
\right.
$$
We notice that
$$
\left.
\begin{array}{l}
Res_1(\omega-\omega_h+w) \leq |\|Res_1|\|_{-1,h} |\|(\psi,\omega-\omega_h+w)|\|_{1,h} \; \forall \; \psi \in H_0^1(\Omega)^2,\\
\\
Res_2(\phi-\phi_h+\beta) \leq
\|Res_2\|_{-1} |\phi-\phi_h+\beta|_{1}.
\end{array}
\right.
$$
 Introducing now the parameters $\nu_1>0$ and $\nu_2>0$ and using two times Young's inequality lead to
$$
\left.
\begin{array}{ll}
(e_h^{rot})^2 &\leq \displaystyle \nu_1 A^2(\varepsilon) |\|Res_1|\|^2_{-1,h}+\frac{1}{\nu_1}|\|(\psi,\omega-\omega_h+w)|\|^2_{1,h}\\[10pt]
&\displaystyle \hspace{5mm} +\nu_2 A^2(\varepsilon) \|Res_2\|^2_{-1}+\frac{1}{\nu_2}|\phi-\phi_h+\beta|^2_{1}\\[10pt]
& \displaystyle \hspace{5mm} -A(\varepsilon)\|\phi-\phi_h+\beta\|^2_{\mathcal{C}}+A(\varepsilon)\|\beta\|^2_{\mathcal{C}}\\[10pt]
&\displaystyle \hspace{5mm} +\left(\frac{\frac{2}{\varepsilon}-1}{1-2\varepsilon}\right)\|\phi_h-\textbf{R}_h\phi_h\|^2 +2\|Res_2\|^2_{-1} +2\|rot(\phi_h-\textbf{R}_h \phi_h)\|^2\\[10pt]
&\displaystyle \hspace{5mm} -\sum_{T \in \mathcal{T}_h}
\left(\frac{\lambda
A(\varepsilon)}{t^2+h_T^2}\right)\|\nabla
(\omega-\omega_h+w)-(\phi-\phi_h+\beta)\|^2_{T}.
\end{array}
\right.
$$
Finally, choosing $\psi=\phi-\phi_h+\beta$, we get
$$
|\|(\psi,\omega-\omega_h+w)|\|^2_{1,h}=\|\nabla(\phi-\phi_h+\beta)\|^2+\sum_{T \in \mathcal{T}_h}\frac{1}{t^2+h_T^2}\|\nabla
(\omega-\omega_h+w)-(\phi-\phi_h+\beta)\|^2_{T},
$$
and so (\ref{relavecA}) holds.
\end{proof}
\bc \label{coro44} Let us assume that $t \leq \sqrt{3 \lambda
c_F^2/\mu}$, and let us define :
\begin{equation}
\label{valeurzeta}
\zeta=\max \left\{ \frac{1}{\mu}, \frac{1}{2 \, \lambda}\right\}.
\end{equation}
Then,
\begin{eqnarray*}
&&(e_h^{rot})^2 \displaystyle \leq 2 \zeta
\left(\frac{3}{\mu}+\frac{c^2_F}{\mu}(3+2\sqrt{3})+4 (\mu + \tilde{\lambda})\right)|\|Res_1|\|^2_{-1,h}\\[10pt]
&&\displaystyle \hspace{5mm}+\left(2 \zeta \left(\frac{3}{\mu}+\frac{c^2_F}{\mu}(3+2\sqrt{3})+4 (\mu + \tilde{\lambda})\right)+2\right)\|Res_2\|^2_{-1}\\[10pt]
&&\displaystyle \hspace{5mm} +\max \left(7+4\sqrt{3} \; ;  \;
2+\left(\frac{3}{\mu}+\frac{c^2_F}{\mu}(3+2\sqrt{3})+4 (\mu +
\tilde{\lambda})\right)2(\mu + \tilde{\lambda})c_{R}^{2}\right)
\|\phi_h-\textbf{R}_h\phi_h\|^2_{H(rot,\Omega)}.
\end{eqnarray*}
\ec
\begin{proof}
Assuming $1-2\varepsilon
>0$, the parameters $\nu_1$ and $\nu_2$ arising in the values of
$A_3$ and $A_5^T$ in (\ref{relavecA}) are first chosen such that
$A_3 \leq 0$ and $A_5^T \geq 0 \; \forall T \in \mathcal{T}_h$. Namely
we take $\nu_1=\nu_2=2 \, \zeta/A(\varepsilon)$. Consequently, we
obtain
\begin{equation}\label{relavecAbis}
\left.
\begin{array}{ll}
(e_h^{rot})^2&\leq \tilde{A}_1 |\|Res_1|\|^2_{-1,h}+\tilde{A}_2 \|Res_2\|^2_{-1}+\tilde{A}_4 \|\phi_h-\textbf{R}_h\phi_h\|^2_{H(rot,\Omega)},
\end{array}
\right.
\end{equation}
with
$$
\left|
\begin{array}{l}
\displaystyle \tilde{A}_1=2 \zeta A(\varepsilon);\\[10pt]
\displaystyle \tilde{A}_2=2 \zeta A(\varepsilon)+2;\\[10pt]
\displaystyle
\tilde{A}_4=\max \left( \frac{\frac{2}{\varepsilon}-1}{1-2\varepsilon} \; ; \; 2+2A(\varepsilon)
(\mu + \tilde{\lambda})c_R^{2} \right).
\end{array}
\right.
$$
Now, in order to provide a result as sharp as possible, it remains
to choose appropriately the parameter $\varepsilon$ to make
the coefficients  $\tilde{A}_1$, $\tilde{A}_2$ and
  $\tilde{A}_4$ arising in (\ref{relavecAbis}) as small as possible. Since we always have $1 \leq  3/\mu+4 (\mu + \tilde{\lambda})$, the assumption $t \leq \sqrt{3 \lambda c_F^2/\mu}$ leads to
 $$
A(\varepsilon)=\displaystyle \frac{3}{\mu}+c_F^2\frac{\displaystyle \frac{1}{\varepsilon}+\varepsilon-1}{\mu(1-2\varepsilon)}+4 (\mu + \tilde{\lambda}).
$$
 At this stage we remark that the two functions $A(\varepsilon)$
 as well as $\frac{\frac{2}{\varepsilon}-1}{1-2\varepsilon} $ reach their minimum value for the same
 value of the  argument $\varepsilon$, namely for $\varepsilon=2-\sqrt{3}$.
So, by a simple calculation,  corollary \ref{coro44} holds.
\end{proof} \\
Now, it remains to bound each of the two residuals. \bl
\label{lemm45} Let $N\in \N^*$ be such that $\displaystyle \max_{T
\in \mathcal{T}_h}Y(T) \leq N$, with $Y(T)=\#\{T'\in \mathcal{T}_h\;
| \; T' \subset \omega_T\}$ and $\omega_T=\{K \in \mathcal{T}_h | K \cap T \neq \emptyset \}$
is the patch of elements surrounding $T$ (consequence of the mesh regularity).
Then there exists  $\kappa_2>0$ only depending on the mesh
regularity such that
\begin{equation}
\label{evalres1}
\displaystyle |\|Res_1|\|^2_{-1,h} \leq 2\,N\, \kappa_2^2\sum_{T
\in \mathcal{T}_h}(t^2+h_T^2)\|\gamma_h-y^*\|^2_{T}+ osc^2(g),
\end{equation}
where $osc(g)$ corresponds to an oscillating term.
\el
\begin{proof}
For any $v \in H_0^1(\Omega)$, let us consider $v_h=Jv$ where $J
:H_0^1(\Omega) \rightarrow W_h$ is defined such that (see, for
example \cite{Clement}, known as the Cl\'ement operator)
\begin{equation}
\label{defk1}
\exists \; \kappa_1 >0 \, ; \, \forall \; T \in \mathcal{T}_h,   \|\nabla v_h\|_{T} \leq \kappa_1 \|\nabla v\|_{\omega_T}.
\end{equation}
Moreover, it can be shown\ \cite{CH08} that there exists $\kappa_2>0$ and
$\kappa_3>0$ such that for all $T \in \mathcal{T}_h$  and for any
$\psi \in H^1_0(\Omega)^2$,
$$
\|\nabla (v-v_h)\|_{T}\leq \kappa_2\left(\|\nabla
v-\psi\|_{\omega_T}+h_T\|\nabla \psi\|_{\omega_T}\right),
$$
$$
h_T^{-1} \|v-v_h\|_{T} \leq \kappa_3\left(
\|\nabla v-\psi\|_{\omega_T}+ h_T \|\nabla
\psi\|_{\omega_T}\right).
$$

Then for all $v \in H_0^1(\Omega)$, we get
$$
\left.
\begin{array}{ll}
Res_1(v)
&=\displaystyle Res_1(v-v_h)\\[8pt]
&=\displaystyle (g,v-v_h)-(\gamma_h,\nabla (v-v_h))\\[8pt]
&=\displaystyle (g+div y^*,v-v_h)-(\gamma_h-y^*,\nabla (v-v_h))\\[8pt]
&=\displaystyle \sum_{T \in \mathcal{T}_h} \left((g+div y^*,v-v_h)_{T}-(\gamma_h-y^*,\nabla (v-v_h))_{T}\right)\\[10pt]
&\leq \displaystyle \sum_{T \in \mathcal{T}_h}h_T\sqrt{t^2+h_T^2}\|g+div y^*\|_{T}\times \frac{h_T^{-1}}{\sqrt{t^2+h_T^2}}\|v-v_h\|_{T}\\[10pt]
&\displaystyle \hspace{5mm} +\sum_{T \in \mathcal{T}_h}\sqrt{t^2+h_T^2}\|\gamma_h-y^*\|_{T}\times \frac{1}{\sqrt{t^2+h_T^2}}\|\nabla(v-v_h)\|_{T}.
\end{array}
\right.
$$

So, we can write
$$
\left.
\begin{array}{ll}
\hspace{-5mm}Res_1(v)&\displaystyle \leq \sum_{T \in
\mathcal{T}_h}h_T\sqrt{t^2+h_T^2}\|g+div y^*\|_{T}\\[8pt]
&\displaystyle \hspace{3cm}\times \frac{\kappa_3}{\sqrt{t^2+h_T^2}}\left(\|\nabla
v-\psi\|_{\omega_T}+h_T \|\nabla \psi\|_{\omega_T}\right)\\[10pt]
&\displaystyle \hspace{5mm} +\sum_{T \in
\mathcal{T}_h}\sqrt{t^2+h_T^2}\|\gamma_h-y^*\|_{T}\\[8pt]
&\displaystyle \hspace{3cm}\times
\frac{\kappa_2}{\sqrt{t^2+h_T^2}}\left(\|\nabla
v-\psi\|_{\omega_T}+h_T
\|\nabla \psi\|_{\omega_T}\right)\\[10pt]
&\displaystyle \leq \left(\sum_{T \in
\mathcal{T}_h}\kappa_3^2 h_T^2(t^2+h_T^2)\|g+div
y^*\|^2_{T}\right. \\[8pt]
&\displaystyle \hspace{3cm}+\left. \sum_{T \in
\mathcal{T}_h}\kappa_2^2(t^2+h_T^2)\|\gamma_h-y^*\|^2_{T}\right)^{1/2}\\[8pt]
&\displaystyle \hspace{5mm}\left(2\underbrace{\sum_{T \in
\mathcal{T}_h}\left(\frac{1}{t^2+h_T^2}\|\nabla
v-\psi\|^2_{\omega_T}+\frac{h_T^2}{t^2+h_T^2}\|\nabla
\psi\|^2_{\omega_T}\right)}_{=S}\right)^{1/2}.
\end{array}
\right.
$$
Now  recalling  that $\displaystyle \max_{T \in \mathcal{T}_h}Y(T)
\leq N$ we have
$$
\left.
\begin{array}{ll}
S&\displaystyle \leq N \sum_{T \in
\mathcal{T}_h}\left(\frac{1}{t^2+h_T^2}\|\nabla
v-\psi\|^2_{T}+\underbrace{\frac{h_T^2}{t^2+h_T^2}}_{ \leq
1}\|\nabla \psi\|^2_{T}\right)\\[10pt]
&\displaystyle \leq N \sum_{T \in
\mathcal{T}_h}\left(\frac{1}{t^2+h_T^2}\|\nabla
v-\psi\|^2_{T}+\|\nabla \psi\|^2_{T}\right)\\[10pt]
&\displaystyle \leq N \left(\|\nabla \psi\|^2_{\Omega}+\sum_{T
\in \mathcal{T}_h}\frac{1}{t^2+h_T^2}\|\nabla
v-\psi\|^2_{T}\right)\\[10pt]
S&\displaystyle \leq N |\|(\psi,v)|\|^2_{1,h}.
\end{array}
\right.
$$
So we get
$$
\left.
\begin{array}{ll}
Res_1(v)&\displaystyle \leq \left(\kappa_3^2 \sum_{T
\in \mathcal{T}_h}h_T^2(t^2+h_T^2)\|g+div
y^*\|^2_{T} \right.\\[12pt]
&\displaystyle \hspace{10mm}+\left. \kappa_2^2\sum_{T \in
\mathcal{T}_h}(t^2+h_T^2)\|\gamma_h-y^*\|^2_{T}\right)^{1/2} \times
\sqrt{2N} |\|(\psi,v)|\|_{1,h}.
\end{array}
\right.
$$
Consequently
$$
\left.
\begin{array}{ll}
|\|Res_1|\|^2_{-1,h}&\displaystyle \leq 2\, N \,  \left(\kappa_3^2 \, \sum_{T \in
\mathcal{T}_h} h_T^2(t^2+h_T^2)\|g+div
y^*\|^2_{T} \right.\\[10pt]
&\displaystyle \hspace{2cm}\left. +\kappa_2^2 \, \sum_{T \in
\mathcal{T}_h}(t^2+h_T^2)\|\gamma_h-y^*\|^2_{T} \right).
\end{array}
\right.
$$
Since $div\, y^*=-\Pi_h \, g$, we get $\|g+div y^*\|^2_{T} \leq C\,h_T^2
\|g\|^2_{\omega_T}$ and  (\ref{evalres1}) holds.
\end{proof}

\bl \label{lemm46} For $\Psi \in H_0^1(\Omega)^2$, we have
\begin{equation}
\label{res2val}
\begin{array}{lcl}
Res_2(\psi) &\leq& \|\mathcal{C}^{-1/2} (x^*-\mathcal{C}
\varepsilon (\phi_h))\| \;  \|\psi\|_{\mathcal{C}}.
\end{array}
\end{equation}
\el
\begin{proof}
Using standard Green formula, we easily obtain
$$
Res_2(\psi)=\int_{\Omega} (x^*-\mathcal{C} \varepsilon (\phi_h))
: \varepsilon (\psi) \; dx +\int_{\Omega}(\gamma_h+div \, x^*) \; \psi
\, dx,
$$
Since $\mathcal{C}$ is a symmetric positive definite operator, we
can define $\mathcal{C}^{1/2}$ and $\mathcal{C}^{-1/2}$ such that
$\mathcal{C}^{1/2} \, \circ \, \mathcal{C}^{1/2}=\mathcal{C}$ and
$\mathcal{C}^{1/2}\circ \mathcal{C}^{-1/2} = \mathcal{I}$. Then the
definition of $x^*$ directly yields
$$
Res_2(\psi)=\int_{\Omega} \mathcal{C}^{-1/2}(x^*-\mathcal{C}
\varepsilon (\phi_h)) \,:\,
\mathcal{C}^{1/2}\varepsilon (\psi) \; dx,
$$
and the Cauchy-Schwarz inequality finally leads to (\ref{res2val}).
\end{proof}
\bt[Reliability of the estimator] Under the assumption of corollary
\ref{coro44}, we have \label{maintheo111}
$$
\left.
\begin{array}{l}
(e_h^{rot})^2 \leq 4 \, \zeta\, N \, \kappa^2_2
\left(\displaystyle \frac{3}{\mu}+\frac{c^2_F}{\mu}(3+2\sqrt{3})+4 (\mu + \tilde{\lambda})\right)\displaystyle \sum_{T \in
\mathcal{T}_h}(t^2+h_T^2)\|\gamma_h-y^*\|^2_{T}\\[10pt]
\displaystyle \hspace{5mm}+4 \left(\zeta \,
\left(\frac{3}{\mu}+\frac{c^2_F}{\mu}(3+2\sqrt{3})+4 (\mu + \tilde{\lambda})\right)+1 \right) (\mu + \tilde{\lambda})
\|\mathcal{C}^{-1/2} (x^*-\mathcal{C} \varepsilon
(\phi_h))\|^2\\[10pt]
\displaystyle \hspace{5mm} +\max
 \left((7+4\sqrt{3}) ;2+ 2 \left(\frac{3}{\mu}+\frac{c^2_F}{\mu}(3+2\sqrt{3})+4 (\mu + \tilde{\lambda})\right) (\mu + \tilde{\lambda}) c_R^2 \right) \|\phi_h-\textbf{R}_h\phi_h\|^2_{H(rot,\Omega)}\\[10pt]
\displaystyle \hspace{5cm} +osc^2(g).
\end{array}
\right.
$$
\et
\begin{proof}
The theorem is a direct consequence of corollary \ref{coro44}, lemma
\ref{lemm45} and lemma \ref{lemm46}.
\end{proof}

\br{\rm In theorem \ref{maintheo111}, all constants are explicitly
given. Indeed, even if $c_F$ and $c_R$ depend on the domain $\Omega$
whereas $\kappa_2$ and $N$ depend on the  used mesh, they can   be
evaluated or at least bounded, see \cite{CF00} and section \ref{sec-num} below devoted to
the numerical validations. }\er

\br{\rm The assumption  $t \leq \sqrt{3 \lambda c^2_F/\mu}$ needed
in corollary \ref{coro44} is not restrictive
since, in the Reissner-Mindlin model, $t$ is expected to be a
very small parameter, so that this property is naturally obtained. }
\er

\section{Efficiency of the estimator}
\label{efficiency}
In order to prove the efficiency of the estimator, each part of it has now to be bounded by the error $e_h^{rot}$ up to a multiplicative constant. In the following, the notation $a\lesssim b$ and $a\sim b$ means the existence of positive constants $c_1$ and $c_2$, which are
independent of the mesh size, of the plate thickness parameter $t$, of the quantities $a$ and $b$ under
consideration  and of the coefficients of the operators such that
$a\le c_2 \, b$ and $c_1 \, b\le a\le c_2 \, b$, respectively. The constants may in particular depend on the aspect ratio $\sigma$ of the mesh.
\bl
\label{rh-i}
$$
\left.
\begin{array}{ll}
\displaystyle \|(\textbf{R}_h-I)\phi_h\|^2_{H(rot,\Omega)}&\lesssim
\lambda^{-2} t^4
\|\gamma-\gamma_h\|^2_{\Omega}+|\omega-\omega_h|^2_{1}\\[8pt]
&\displaystyle
\hspace{2cm}+|\phi-\phi_h|^2_{1}+\lambda^{-2} t^4
\|rot(\gamma-\gamma_h)\|^2.
\end{array}
\right.
$$
\el
\begin{proof}
Since
$$
\displaystyle (\textbf{R}_h-I)\phi_h = \lambda^{-1} t^2
(\gamma-\gamma_h)-\nabla
(\omega-\omega_h)+(\phi-\phi_h),
$$
we have
$$
\|\textbf{R}_h-I)\phi_h\|
\leq \lambda^{-1} t^2 \|\gamma-\gamma_h\|+|\omega-\omega_h|_{1}+\|\phi-\phi_h\|,
$$
and with the Poincar\'e-Friedrichs inequality, we get
$$
\displaystyle \|(\textbf{R}_h-I)\phi_h\|^2 \lesssim
\lambda^{-2} t^4
\|\gamma-\gamma_h\|^2+|\omega-\omega_h|^2_{1}+|\phi-\phi_h|^2_{1}.
$$
Moreover, we have
$$
\left.
\begin{array}{l}
\|rot(\phi_h-\textbf{R}_h\phi_h)\|^2 \lesssim \lambda^{-2} t^4
\|rot(\gamma-\gamma_h)\|^2+|\phi-\phi_h|^2_{1},
\end{array}
\right.
$$
so that lemma \ref{rh-i} holds.
\end{proof} \\

\bl \label{lemmpart2} There exists a relevant choice of $x^*$ such
that \be \label{majsecondpart}
\|\mathcal{C}^{-1/2}(x^*-\mathcal{C}\varepsilon(\phi_h))\|^2
\lesssim \|\gamma_h-\gamma\|^2_{-1}+|\phi_h-\phi|^2_1. \ee \el
\begin{proof}
First, there exists only one pair $(\phi_h^{*},\phi_h^{**}) \in
H_0^1(\Omega)^2 \times \Theta_h$ solution of
$$
\left\{
\begin{array}{lcll}
a(\phi_h^{*},\psi)&=&-(\gamma_h,\psi) &\; \forall \; \psi \in H_0^1(\Omega)^2,\\[10pt]
a(\phi_h^{**},\psi_h)&=&-(\gamma_h,\psi_h) &\; \forall \; \psi_h \in \Theta_h.
\end{array}
\right.
$$
Then, by Theorem 3.9 of \cite{nicaise-wohlmuth:05} and a relevant
construction
 of $x^*$, for all $T$ in $\mathcal{T}_h$ we have
$$
\displaystyle
\begin{array}{lcl}
\|\mathcal{C}^{-1/2}(x^*-\mathcal{C}\varepsilon(\phi_h^{**}))\|_{T}
&\lesssim& \|\phi_h^*-\phi_h^{**}\|_{\mathcal{C}, \omega_T}.
\end{array}
$$
Because of the mesh regularity, we also get the global estimate
\begin{equation}\label{phi*}
\displaystyle
\|\mathcal{C}^{-1/2}(x^*-\mathcal{C}\varepsilon(\phi_h^{**}))\|
\lesssim \|\phi_h^*-\phi_h^{**}\|_{\mathcal{C}}.
\end{equation}
Clearly
$$
\displaystyle
\mathcal{C}^{-1/2}(x^*-\mathcal{C}\varepsilon(\phi_h))=
\mathcal{C}^{-1/2}(x^*-\mathcal{C}\varepsilon(\phi_h^{**}))+\mathcal{C}^{1/2}\varepsilon(\phi_h^{**}-\phi_h).
$$
By (\ref{phi*}) and the triangular inequality, we  arrive at
\begin{eqnarray}
\|\mathcal{C}^{-1/2}(x^*-\mathcal{C}\varepsilon(\phi_h))\|
&\lesssim
\|\mathcal{C}^{-1/2}(x^*-\mathcal{C}\varepsilon(\phi_h^{**}))\|+\|\phi_h^{**}-\phi_h\|_{\mathcal{C}} \nonumber\\[10pt]
\displaystyle &\lesssim \|\phi_h^*-\phi_h^{**}\|_{\mathcal{C}}+
\|\phi_h^{**}-\phi_h\|_{\mathcal{C}}. \label{rel22}
\end{eqnarray}
Now, it remains to bound each of the two terms of the right-hand side of (\ref{rel22}).
To begin with, let us consider $\psi_h \in \Theta_h$. Thanks to the definition of $\phi_h^{**}$, we get
$$
\left.
\begin{array}{ll}
a(\phi_h-\phi_h^{**},\psi_h)&=(\gamma_h,\psi_h-\textbf{R}_h
\psi_h)\\[8pt]
&=\displaystyle
(\gamma_h,\psi_h)-(\gamma_h,\textbf{R}_h
\psi_h)\\[8pt]
&=\displaystyle
(\gamma_h-\gamma,\psi_h)+a(\phi_h-\phi,\psi_h)\\[8pt]
&\lesssim \displaystyle
(\|\gamma_h-\gamma\|_{-1}+|\phi_h-\phi|_{1})
|\psi_h|_{1}.
\end{array}
\right.
$$
By taking $\psi_h=\phi_h-\phi_h^{**},$ we obtain
\begin{equation} \label{**}
\|\phi_h^{**}-\phi_h\|_{\mathcal{C}} \lesssim \displaystyle
\|\gamma_h-\gamma\|_{-1}+|\phi_h-\phi|_{1}.
\end{equation}
Then, by the triangular inequality, we get
$$
\|\phi_h^*-\phi_h\|_{\mathcal{C}} \leq
\|\phi_h^*-\phi\|_{\mathcal{C}} + \|\phi-\phi_h\|_{\mathcal{C}},
$$
and by the definition of $\phi_h^*$, we have for all $\psi \in H_0^1(\Omega)^2$
$$
a(\phi_h^*-\phi,\psi)=(\gamma-\gamma_h,\psi),
$$
so that
$$
\|\phi_h^*-\phi\|_{\mathcal{C}} \lesssim
\|\gamma-\gamma_h\|_{-1}.
$$
We then obtain
\begin{equation}\label{*}
\|\phi_h^*-\phi_h\|_{\mathcal{C}} \leq
\|\gamma-\gamma_h \|_{-1} +
\|\phi-\phi_h\|_{\mathcal{C}}\lesssim
\|\gamma-\gamma_h\|_{-1} +
|\phi-\phi_h|_1.
\end{equation}
Using (\ref{**}) and (\ref{*}) in (\ref{rel22}), we get (\ref{majsecondpart}).
\end{proof}
\bl \label{troislem} There exists a relevant choice of $y^*$ such
that
\begin{equation}
 \label{troislemeq}
\displaystyle \sum_{T \in \mathcal{T}_h}(t^2+h_T^2)
\|\gamma_h-y^*\|^2_{T} \lesssim
t^2\|\gamma-\gamma_h\|^2
+\|\gamma-\gamma_h\|_{-1}^2+osc^2(g),
\end{equation}
where $osc^2(g)$ is an oscillating term.
\el
\begin{proof}
Because of lemma 3.1 of \cite{DN08}, we have for any  $T \in
\mathcal{T}_h$ the equivalence
$$
\|\gamma_h-y^*\|_{T} \sim \displaystyle  h_T^{1/2} \sum_{E \subset
\partial T}\|(\gamma_h-y^*)\cdot \nu_T\|_{E},
$$
where $\nu_T$ is the outward unit normal vector to $T$. Now we
define $y^*$ as in \cite{DN08}, by noticing that (\ref{pbdiscret})
implies that $$(\gamma_h, \nabla v_h) = (g,v_h) \quad  \forall v_h  \in
W_h,
$$
hence there exist fluxes $g_E\in \P_1(E),$ for all edges $E$ such
that
$$\int_T\gamma_h\cdot \nabla v_h =\int_T g v_h+ \int_{\partial T} g_T v_h \quad \forall v_h  \in \P_1(T),
$$
where $g_T=g_E \nu_E \nu_T$, $\nu_E$ being a fixed normal vector to
$E$. According to the definition of the $BDM_1$ elements there then
exists a unique $y_T^*\in \P_1(T)^2$ such that
\[
y_T^*\cdot \nu_E= g_E \quad \forall E\subset T.
\]
Hence we define $y^*$ such that its restriction to each triangle $T$
is equal to $y_T^*$. According to its definition $y^*$ belongs to
$H_{div}(\Omega)$ and  moreover according to \cite{DN08}, we have
\[
\div y^*=-\Pi_h g.
\]
Then by the use of theorem 6.2 from \cite{AinsworthOden} and the
mesh regularity we get
$$
 \|\gamma_h-y^*\|_{T} \lesssim h_T^{1/2} \sum_{E
\subset
\partial T \setminus \partial \Omega}\|[\gamma_h \cdot \nu_E]_E\|_{E}+\sum_{T' \subset \omega_T} h_{T'} \|div \, \gamma_h +g\|_{T'},
$$
where $[ v ]_E$ denotes the jump of the quantity
$v$ through the edge $E$. Consequently
\begin{equation}
\label{eq36}
\begin{array}{lcl}
\displaystyle \sum_{T \in \mathcal{T}_h}(t^2+h_T^2)
\|\gamma_h-y^*\|^2_{T} &\lesssim& \displaystyle \sum_{T \in \mathcal{T}_h}h_T
(t^2+h_T^2) \sum_{E \subset
\partial T \setminus \partial \Omega}\|[\gamma_h \cdot \nu_E]_E\|^2_{E}\\
&+&\displaystyle \sum_{T \in \mathcal{T}_h}\sum_{T'
\subset \omega_T} h^2_{T'}(t^2+h_T^2) \|div \gamma_h
+g\|^2_{T'} \\
 &\lesssim& \displaystyle \sum_{E \subset
\partial T \setminus \partial \Omega}h_E (t^2+h_E^2)\|[\gamma_h \cdot \nu_E]_E\|^2_{E} \\
&+&\displaystyle \sum_{T \in \mathcal{T}_h} h^2_{T}(t^2+h_T^2) \|div
\gamma_h +g\|^2_{T}.
\end{array}
\end{equation}

\noindent
Using the classical edge bubble functions as well as elementwise inverse estimates, it is proved in \cite{CH08}, section 4.3 that :
\begin{equation}
\label{term1}
\left.
\begin{array}{l}
\displaystyle \sum_{E \in \mathcal{E}(\Omega)\setminus
\partial \Omega} h_E (t^2+h_E^2) \|[\gamma_h \cdot \nu_E]_E\|^2_{E} \lesssim
\sum_{T \in \mathcal{T}_h}
h_T^2(t^2+h_T^2)\|g-\Pi_h g\|^2_{T}\\[8pt]
\displaystyle \hspace{9cm}+
\|\gamma-\gamma_h\|_{-1}^2+t^2\|\gamma-\gamma_h\|^2.
\end{array}
\right.
\end{equation}
\noindent
Moreover, with the use of classical element bubble functions as well as elementwise inverse estimates, it is also  proved in \cite{CH08}, section 4.1 that :
\begin{equation}\label{divgamma}
\left.
\begin{array}{lcl}
\displaystyle \sum_{T \in \mathcal{T}_h}
h_T^2(t^2+h_T^2) \|div \gamma_h+\Pi_hg\|^2_{T} & \lesssim &
t^2\|\gamma-\gamma_h\|^2 + \|\gamma-\gamma_h\|^2_{-1}\\
&+&\displaystyle \sum_{T \in \mathcal{T}_h}h_T^2(h_T^2+t^2)\|g-\Pi_hg\|^2_{T}.
\end{array}
\right.
\end{equation}
Now, from (\ref{eq36})  associated to the standard triangular inequality :
$$
\|div\gamma_h+g\|_{T} \leq
\|div\gamma_h+\Pi_hg\|_{T}+\|g-\Pi_h g\|_{T},
$$
the use of (\ref{term1}) and (\ref{divgamma}) leads to (\ref{troislemeq}).

\end{proof}
\bt[Efficiency of the estimator] There exists a relevant choice of
$x^*$ and of $y^*$ such that \label{effesti2}
$$
\left.
\begin{array}{l}
\displaystyle 4 \, \zeta \, N \, \kappa^2_2
\left(\frac{3}{\mu}+\frac{c^2_F}{\mu}(3+2\sqrt{3})+4 (\mu + \tilde{\lambda})\right)\sum_{T \in
\mathcal{T}_h}(t^2+h_T^2)\|\gamma_h-y^*\|^2_{T}\\[10pt]
\displaystyle \hspace{5mm}+\left(2 \, \zeta \,
(\frac{3}{\mu}+\frac{c^2_F}{\mu}(3+2\sqrt{3})+4 (\mu + \tilde{\lambda}))+2\right)2 (\mu + \tilde{\lambda})
\|\mathcal{C}^{-1/2} (x^*-\mathcal{C} \varepsilon
(\phi_h))\|^2\\[10pt]
\displaystyle \hspace{5mm} + \max
\left( (7+4\sqrt{3});2+(\frac{3}{\mu}+\frac{c^2_F}{\mu}(3+2\sqrt{3})+4 (\mu + \tilde{\lambda}))2 (\mu + \tilde{\lambda}) c_R^2 \right) \|\phi_h-\textbf{R}_h\phi_h\|^2_{H(rot,\Omega)}\\[10pt]
\displaystyle  \lesssim (e_h^{rot})^2+ osc^2(g).
\end{array}
\right.
$$
\et
\begin{proof}
The proof is a direct consequence of lemma \ref{rh-i}, \ref{lemmpart2} and \ref{troislem}.
\end{proof}

\section{Numerical validation}
\label{sec-num} Here we illustrate and validate our theoretical
results by a simple computational example. Let $\Omega$ be the unit
square $]0,1[^2$. We consider the exact solution $(\omega, \phi)$ in
$\Omega$ of the Reissner-Mindlin problem (\ref{S2})-(\ref{S3}) given
by
$$
\displaystyle \phi= \left(
\begin{array}{l}
\displaystyle\frac{1-2x}{x^2(1-x)^2}\\[10pt]
\displaystyle\frac{1-2y}{y^2(1-y)^2}
\end{array}
\right) exp\left(-\frac{1}{x(1-x)}-\frac{1}{y(1-y)}\right),
$$
and
$$
\displaystyle \omega=\left(1-(2\mu+\tilde{\lambda}) \lambda^{-1} t^2 \;(a(x)+a(y))\right)
\;exp\left(-\frac{1}{x(1-x)}-\frac{1}{y(1-y)}\right),
$$
with
$$
a(z)=\displaystyle \frac{6z^4-12z^3+12z^2-6z+1}{z^4(1-z)^4}.
$$
The corresponding scaled transverse loading function $g$ is given by
$$
g=(2\mu+\tilde{\lambda})\;\left(c(x)+c(y)+2\,a(x)\,a(y)\right)\;
exp \left(-\frac{1}{x(1-x)}-\frac{1}{y(1-y)}\right)
$$
with
$$
\begin{array}{lcl}
c(z)&=&\displaystyle \frac{120z^{10}-600z^9+1620z^8-2880z^7+3504z^6-2952z^5}{z^8(1-z)^8} \\
&+&\displaystyle \frac{1708z^4-656z^3+156z^2-20z+1}{z^8(1-z)^8}.
\end{array}
$$
This analytical solution is extended by $0$ on $\partial \Omega$ to
obtain $(\omega, \phi) \in H_0^1(\Omega) \times H_0^1(\Omega)^2$.
Here we take $t=1/1024$, $\lambda=1$, $\mu= 1$ and
$\tilde{\lambda}=1$. The meshes  we use are uniform ones composed of $n^2$ squares, each of
them being cut into 8 triangles as displayed on Figure
\ref{figuremaillage} for $n=4$. The refinement strategy is an
uniform one  so that the   value of the mesh size $h$ between two
consecutive meshes is twice smaller.
\begin{figure}[h]
\centering
\includegraphics[width=4cm]{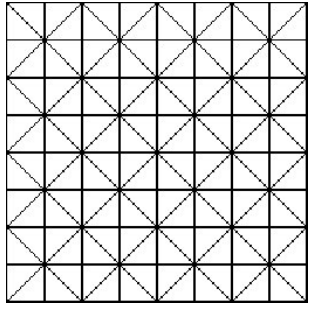}
\caption{Mesh level corresponding to $n=4$ and $h=\sqrt{2}/8$.\label{figuremaillage}}
\end{figure}
In order to validate the reliability of the estimator, we
consider the ''discrete error'' given by
\begin{eqnarray}
 e_{h,dis}^{rot}&=&\sqrt{|\omega - \omega_h|^2_{1} +
|\phi-\phi_h|^2_{1} + \lambda^{-1} t^2
\|\gamma-\gamma_h\|^2 + \lambda^{-2} t^4
\|rot(\gamma-\gamma_h)\|^2 + ||P_h \gamma-\gamma_h||_{-1,h}^2}, \nonumber
\end{eqnarray}
where $P_h \gamma$ stands for the  piecewise $\Poly_1$-discontinuous
interpolation of $\gamma$ on the mesh $\mathcal{T}_h$.  This
discrete error is defined by approximating the $H^{-1}(\Omega)$ norm
of $\gamma-\gamma_h$ arising in $e_{h}^{rot}$ (see (\ref{ehrot2})) by its discrete
locally computable version defined by
\begin{equation}
||P_h \gamma-\gamma_h||_{-1,h}^2 = \displaystyle \sup_{v_h \in W_h } \frac{|(P_h \gamma-\gamma_h,v_h)|^2}{|v_h|^2_1}.
\end{equation}
The computation of $||P_h \gamma-\gamma_h||_{-1,h}^2$ is now an easy
task  and simply corresponds to the determination of the largest
eigenvalue of a classical generalized finite dimensional eigenvalue
problem. In order to validate the reliability of the estimator
according to theorem \ref{maintheo111}, the error estimator is
defined by
$$
\left.
\begin{array}{l}
(\eta_h)^2 =  4 \, \zeta\, N \, \kappa^2_2
\left(\displaystyle \frac{3}{\mu}+\frac{c^2_F}{\mu}(3+2\sqrt{3})+4 (\mu + \tilde{\lambda})\right)\displaystyle \sum_{T \in
\mathcal{T}_h}(t^2+h_T^2)\|\gamma_h-y^*\|^2_{T}\\[10pt]
\displaystyle \hspace{5mm}+4 \left(\zeta \,
\left(\frac{3}{\mu}+\frac{c^2_F}{\mu}(3+2\sqrt{3})+4 (\mu + \tilde{\lambda})\right)+1 \right) (\mu + \tilde{\lambda})
\|\mathcal{C}^{-1/2} (x^*-\mathcal{C} \varepsilon
(\phi_h))\|^2\\[10pt]
\displaystyle \hspace{5mm} +\max
 \left((7+4\sqrt{3}) ;2+ 2 \left(\frac{3}{\mu}+\frac{c^2_F}{\mu}(3+2\sqrt{3})+4 (\mu + \tilde{\lambda})\right) (\mu + \tilde{\lambda}) c_R^2 \right) \|\phi_h-\textbf{R}_h\phi_h\|^2_{H(rot,\Omega)},\\[10pt]
\end{array}
\right.
$$
and we plot on Figure  \ref{figure:fiabilite-esti2} the evolution of
the computed effectivity index $\eta_h/e_{h,dis}^{rot}$ versus $h$.
Here, the values of $x^\star$ as well as $y^\star$ are respectively
computed in the same manner as in  \cite{DN08} and
\cite{nicaise-wohlmuth:05}, in order to obtain relevant choices 
as required by theorem \ref{effesti2} to ensure the efficiency of the estimator.
Practically, some fluxes $g_E$ through the edges $E$ of
each triangle of the mesh are needed, and have to be computed by
solving  local linear problems. In fact, in our tests, these values
are explicitely defined. For $y^*$, we use $g_E =  \{\{{\gamma_h \cdot \nu_E}\}\}$, where $ \{\{{\gamma_h \cdot \nu_E}\}\} $ 
denotes the averaged value on the triangles on each side of $E$ of $\gamma_h \cdot \nu_E$ evaluated at the middle of $E$.
For $x^*$, we use $g_E=\sum_{x \in \mathcal{N}(T)} \{ \{ {{\mathcal{C}\varepsilon (\phi_h)\} \} (x)
 \nu_E}} \lambda_x$.  Here, $\{ \{ \mathcal{C}\varepsilon (\phi_h)\} \} (x)$ is the averaged value over the triangles surrouding the node $x$ of the piecewise constant function on each triangle 
$\mathcal{C}\varepsilon (\phi_h)$, and $\lambda_x$ stands for the classical $\Poly_1$-Lagrange basis function associated with the node $x$.
Moreoever, for the construction of $x^*$, the
Argyris basis functions have to be used  (see section 4 of \cite{nicaise-wohlmuth:05} as well as \cite{DS06} for the practical
implementation). \\
 \\
\noindent
From (\ref{valeurzeta}) we have $\zeta=1$.
The Poincar\'e-Friedrichs constant $c_F$ is here equal to
$1/(\sqrt{2} \pi)$ since $\Omega$ is the unit square. Because of the
kind of meshes used (see Figure \ref{figuremaillage}), we have $N=8$
and $\kappa_2=1+\displaystyle \frac{12}{\sqrt{2}\pi}$ (see annex
\ref{estikappa2}). Finally, it can be proved \cite{HorganPayne:83}
that on the unit square, $c_R\leq 2\sqrt{\frac{1}{2-\sqrt{2}}}$,
hence below we take this upper bound for $c_R$ (while it is
conjectured that $c_R=\sqrt{\frac{2 \displaystyle \pi}{\displaystyle
\pi-2}}$, see \cite{Dau08}).
\begin{figure}[htp]
\centering
\includegraphics[width=10cm]{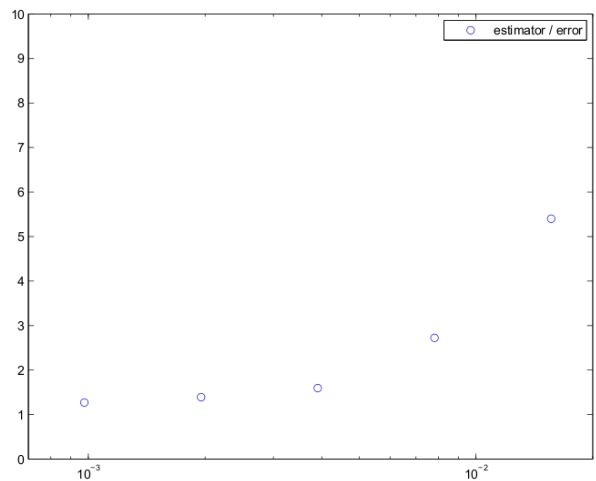}
\vskip-2cm
\caption{$\eta_h/e_{h,dis}^{rot} $ versus $h$. \label{figure:fiabilite-esti2}}
\end{figure}
\noindent As expected by the theory, it can be observed that the
computed effectivity index is larger than one. Moreover, it  converges
towards a constant close to one when $h$ goes towards zero, so that the proposed
estimator is asymptotically exact.
\clearpage
\section{Annexes}
\subsection{Proof of (\ref{Korn})}
\label{annexe1} Let us consider $v \in C_c^{\infty}(\Omega)^2$. Two integrations by parts yield :
$$
\left.
\begin{array}{ll}
\displaystyle 2 \int_{\Omega} |\varepsilon (v)|^2 dx &\displaystyle
= \int_{\Omega} |\nabla v|^2 dx + \int_{\Omega} \nabla v  (\nabla
v)^T dx \\[8pt]
&\displaystyle = \int_{\Omega} |\nabla v|^2 dx + \int_{\Omega}
|div v|^2  dx \\[8pt]
&\displaystyle  \geq \int_{\Omega} |\nabla v|^2 dx.
\end{array}
\right.
$$
Hence by a density argument we obtain
$$
\|\nabla v\|  \leq \sqrt{2} \; \|\varepsilon (v)\| \quad \forall \;
v \in H^1_0(\Omega)^2.
$$
Then, we recall
$$
\displaystyle \mathcal{C} \varepsilon (\phi) \displaystyle= 2 \mu
\varepsilon (\phi) +
\tilde{\lambda} Tr (\varepsilon (\phi)) \mathcal{I},
$$
so that
$$
\left.
\begin{array}{rl}
\displaystyle \|\phi\|^2_{\mathcal{C}}&\displaystyle =
\int_{\Omega} \mathcal{C}
\varepsilon (\phi)  \varepsilon (\phi) dx \\[8pt]
&\displaystyle=2 \mu \int_{\Omega} \varepsilon(\phi)  \varepsilon
(\phi) dx + \tilde{\lambda} \int_{\Omega} Tr (\varepsilon (\phi))
 \mathcal{I}   \varepsilon(\phi) dx \\[8pt]
&\displaystyle=2 \mu \int_{\Omega} |\varepsilon(\phi)|^2 dx +
\tilde{\lambda}
\int_{\Omega} (Tr \varepsilon (\phi))^2 dx \\[10pt]
&\displaystyle \geq  \mu  \|\nabla
\phi\|^2.
\end{array}
\right.
$$
This proves (\ref{Korn}).

\subsection{Evaluation of $\kappa_2$ for the triangulation of section 6}
\label{estikappa2} With the definitions given above, let us consider
$z$ an affine function on $\omega_T$, so that $Jz=z$ on $T$. With
$v$ and $v_h$ defined in the proof of lemma \ref{lemm45} and the
triangular inequality, we get
\begin{eqnarray*}
\|\nabla(v-v_h)\|_{T}& \leq & \displaystyle \|\nabla (v-z)\|_{T}+\|\nabla J(v-z)\|_{T}\\
\end{eqnarray*}
From (\ref{defk1}), we get
\begin{eqnarray*}
\|\nabla(v-v_h)\|_{T}& \leq & \displaystyle(1+\kappa_1)\|\nabla (v-z)\|_{\omega_T}.\
\end{eqnarray*}

\noindent Defining $A=\nabla z$ and considering $\psi \in
H^1_0(\Omega)^2$, we have
$$
\displaystyle \|\nabla (v-v_h)\|_{T}\leq (1+\kappa_1) \|\nabla
v-A\|_{\omega_T} \leq (1+\kappa_1) \left( \|\nabla
v-\psi\|_{\omega_T}+\|\psi-A\|_{\omega_T} \right).
$$
Now, $z$ is chosen such that
$$
A=\frac{1}{|\omega_T|} \int_{\omega_T} \psi dx.
$$
By Poincar\'e inequality, there exists $C_{\omega_T} >0$, depending
on the patch $\omega_T$, such that
$$
\displaystyle
\|\psi-A\|_{\omega_T} \leq  C_{\omega_T} \, h_T \, \|\nabla
\psi\|_{\omega_T}\; \forall \; \psi \in H^1_0(\Omega)^2.
$$
So,
\begin{eqnarray}
\displaystyle \|\nabla (v-v_h)\|_{T}&\leq & \displaystyle(1+\kappa_1) \|\nabla v-\psi\|_{\omega_T}+
(1+\kappa_1)  \, C_{\omega_T} \, h_T \|\nabla
\psi\|_{\omega_T} \nonumber \\[10pt]
&\leq& \displaystyle \underbrace{(1+\kappa_1) max \{1;C_{\omega_T}\}}_{\displaystyle =\kappa_2} \left(
\|\nabla v-\psi\|_{\omega_T}+ h_T \|\nabla
\psi\|_{\omega_T}\right). \label{valkappa2}
\end{eqnarray}
Now, it remains to evaluate $\kappa_1$ as well as $C_{\omega_T}$.
Let $\eta_z$ be the nodal basis associated to $W_h$. We have
$$
J\,v\,=\,\sum_{z \in \mathcal{N}}v_z\,\eta_z\, , \, \forall
\,v\,\in\, H^1_0(\Omega),
$$
from what we deduce
$$
\nabla\,J\,v\,=\,\sum_{z \in \mathcal{N}}(v_z-v)\,\nabla\,\eta_z\, ,
\, \forall \,v\,\in\, H^1_0(\Omega).
$$
Let us define $\mathcal{N}_T=\mathcal{N} \cap T$. We have
\begin{eqnarray*}
\|\nabla\,J\,v\|_T&=&\|\sum_{z
\in \mathcal{N}_T}(v_z-v)\,\nabla\,\eta_z\|_T\\
&\leq&\sum_{z \in
\mathcal{N}_T} \|v_z-v\|_T\,\|\nabla\,\eta_z\|_T\\
&\leq&\sum_{z \in
\mathcal{N}_T}\|v_z-v\|_{\omega_z}\,\|\nabla\,\eta_z\|_T
\end{eqnarray*}
But
$$
\|\nabla\,\eta_z\|_T \leq \rho_T^{-1},
$$
and from \cite[(5.12)]{CF00}, we get
$$
\|v_z-v\|_{\omega_z}  \leq c(\omega_z,2)\|\nabla\,v\|_{\omega_z}.
$$
With the triangulation involved, we have
$$
c(\omega_z,2) \leq \frac{\sqrt2 \, h_T}{\pi},
$$
and
\begin{eqnarray*}
\|\nabla\,J\,v\|_T&\leq &\frac{3 \sqrt2}{\pi} \, \frac{h_T}{\rho_T}\,  \|\nabla \, v\|_{\omega_z},
\end{eqnarray*}
so that
$$
\kappa_1 \leq \frac{3 \sqrt2}{\pi} \, \frac{h_T}{\rho_T}.
$$
For the  involved triangulation  $h_T/\rho_T=2$  and hence
\begin{equation}
\label{valkappa1}
\kappa_1 \leq \displaystyle \frac{12}{\sqrt2 \, \pi}.
\end{equation}
Since from \cite{CF00}, we have $C_{\omega_T}=\displaystyle
\frac{3}{\pi}$, (\ref{valkappa2}) and   (\ref{valkappa1}) lead to
$$
\kappa_2 \leq 1+\frac{12}{\sqrt2 \, \pi}.
$$


\addcontentsline{toc}{section}{Bibliographie}
\begin{thebibliography}{99}

\bibitem{AinsworthOden}
M.~Ainsworth and J.~Oden.
\newblock {\em A posteriori error estimation in finite element analysis}.
\newblock John Wiley and Sons, 2000.

\bibitem{BS01}
I.~Babu\v{s}ka and T.~Strouboulis.
\newblock {\em The finite element methods and its reliability}.
\newblock Clarendon Press, Oxford.

\bibitem{BBF89}
K.J. Bathe, F. Brezzi and M. Fortin. \emph{Mixed-interpolated elements for Reissner-Mindlin plates}. Int. J. Num. Meths. Engrg, Volume 28, pp 1787-1801, 1989.

\bibitem{BD85}
K.J. Bathe and E. Dvorkin. \emph{A four-node plate bending element based on Mindlin-Reissner plate theory and a mixed interpolation}. Int. J. Num. Meths. Engrg., 21, pp 367-383, 1985.

\bibitem{Beiroetco:08}
L.~Beir{\~a}o~da Veiga, C.~Chinosi, C.~Lovadina, and R.~Stenberg.
\newblock \emph{A-priori and a-posteriori error analysis for a family of
  {R}eissner-{M}indlin plate elements}.
\newblock {BIT}, 48(2):189--213, 2008.

\bibitem{braess:06}
D.~Braess and J.~Sch{\"o}berl.
\newblock \emph{Equilibrated residual error estimator for edge elements}.
\newblock {Math. Comp.}, 77(262):651--672, 2008.

\bibitem{brenner:94}
S.~C. Brenner and L.~R. Scott.
\newblock {\em The mathematical theory of finite element methods}.
\newblock Springer, New York, 1994.

\bibitem{brezzifortin}
F.~Brezzi and M.~Fortin.
\newblock {\em Mixed and hybrid finite element methods}.
\newblock Springer-Verlag, New-York, 1991.

\bibitem{Carstensen:02}
C.~Carstensen.
\newblock \emph{Residual-based a posteriori error estimate for a nonconforming
  {R}eissner-{M}indlin plate finite element}.
\newblock {SIAM J. Numer. Anal.}, 39(6):2034--2044 (electronic), 2002.

\bibitem{CF00} C. Cartensen and S.A Funken.
\emph{Constants in Cl\'ement-interpolation error and residual based
a posteriori error estimates in Finite Element Methods}. East-West
J. Numer. Math., Vol 8, No 3. pp 153-175, 2000.

\bibitem{CH08}
C. Cartensen and Jun Hu. \emph{A posteriori error analysis for
conforming MITC elements for Reissner-Mindlin plates}. Mathematics
of computation, 77,  262, pp 611-632, 2008.

\bibitem{CarstensenSchoberl:06}
C.~Carstensen and J.~Sch{\"o}berl.
\newblock \emph{Residual-based a posteriori error estimate for a mixed {R}ei\ss
  ner-{M}indlin plate finite element method}.
\newblock {Numer. Math.}, 103(2):225--250, 2006.

\bibitem{ciarlet:78}
P.~G. Ciarlet.
\newblock {\em The finite element method for elliptic problems}.
\newblock North-Holland, Amsterdam, 1978.

\bibitem{Clement}
Ph. Clement. \emph{Approximation by finite element functions using local regularization}. RAIRO R2, pp 77-84, 1975.

\bibitem{DN08}
S. Cochez-Dhondt and S. Nicaise. \emph{A posteriori error estimators
based on equilibrated fluxes}. {Computational Methods in Applied
Mathematics}, Vol 10,       pp 49-68, 2010.

\bibitem{Dau08}
http://anum-maths.univ-rennes1.fr/JourneesEquipe/04JEAN/Transparents/Dauge.pdf.

\bibitem{DS06} V. Dom\'inguez and F.-J. Sayas. \emph{A
simple Matlab implementation of the Argyris element}. ACM
Trans. Math. Software 35, no. 2, Art. 16, 11 pp, 2009.

\bibitem{DHHLR03}
R. Duran, E. Hernandez, L. Hervella-Nieto, E. Liberman and R. Rodriguez. \emph{Error estimates for lower-order isoparametric quadrilateral finite elements for plates}. SIAM J. Numer. Anal., 41, pp 1751-1772, 2003.

\bibitem{DL92}
Ricardo Dur\'an and Elsa Liberman. \emph{On mixed finite element
methods for the Reissner-Mindlin plate model}. Mathematics of Computation, 58, 198, pp 561-573, 1992.

\bibitem{ern_nic:07}
A.~Ern, S.~Nicaise, and M.~Vohral{\'{\i}}k.
\newblock \emph{An accurate {$\bold H({\rm div})$} flux reconstruction for
  discontinuous {G}alerkin approximations of elliptic problems}.
\newblock {C. R. Math. Acad. Sci. Paris}, 345(12):709--712, 2007.

\bibitem{Repin06}
M.E. Frolov, P. Neittaanm\"aki and S.I. Repin. \emph{Guaranteed
functional error estimates for the Reissner-Mindlin plate problem}.
Journal of Mathematical Sciences, 132, 4, pp 553-561, 2006.

\bibitem{GR86}
V.~Girault and P.-A. Raviart.
\newblock {\em Finite elements methods for Navier-Stokes equations, Theory and
  Algorithms}.
\newblock Springer Series in Computational Mathematics, Berlin, 1986.

\bibitem{HorganPayne:83}
C.~O. Horgan and L.~E. Payne.
\newblock \emph{On inequalities of {K}orn, {F}riedrichs and {B}abu\v ska-{A}ziz}.
\newblock {Arch. Rational Mech. Anal.}, 82(2):165--179, 1983.

\bibitem{HuHuang:10}
J. Hu and Y. Huang.
\newblock \emph{A posteriori error analysis of finite element methods for Reissner-Mindlin plates}.
SIAM J. Numer. Anal., 47 (6):4446-4472, 2010.

\bibitem{LL83}
P.~Ladev\`eze and D.~Leguillon.
\newblock \emph{Error estimate procedure in the finite element method and
  applications}.
\newblock {SIAM J. Numer. Anal.}, 20:485--509, 1983.

\bibitem{Liberman:01}
E.~Liberman.
\newblock \emph{A posteriori error estimator for a mixed finite element method for
  {R}eissner-{M}indlin plate}.
\newblock {Math. Comp.}, 70(236):1383--1396 (electronic), 2001.

\bibitem{LovadinaStenberg:05}
C.~Lovadina and R.~Stenberg.
\newblock \emph{A posteriori error analysis of the linked interpolation technique for
  plate bending problems}.
\newblock {SIAM J. Numer. Anal.}, 43(5):2227--2249 (electronic), 2005.

\bibitem{LW04}
R.~Luce and B.~Wohlmuth.
\newblock \emph{A local a posteriori error estimator based on equilibrated fluxes}.
\newblock {SIAM J. Numer. Anal.}, 42:1394--1414, 2004.

\bibitem{neittaanmaki:04}
P.~Neittaanma\"aki and S.~Repin.
\newblock {\em Reliable methods for computer simulation: error control and a
  posteriori error estimates.}, volume~33 of {Studies in Mathematics and
  its applications}.
\newblock Elsevier, Amsterdam, 2004.

\bibitem{nicaise-wohlmuth:05}
S.~Nicaise, K.~Witowski, and B.~I. Wohlmuth.
\newblock \emph{An a posteriori error estimator for the {L}am\'e equation based on
  equilibrated fluxes}.
\newblock {IMA J. Numer. Anal.}, 28(2):331--353, 2008.

\bibitem{SS97}
R. Stenberg and M. Suri. \emph{An hp error analysis of MITC plate elements}. SIAM J. Numer. Anal. 34, pp 544-568, 1997.

\bibitem{Ver96}
R.~Verfurth.
\newblock {\em A review of a posteriori error estimation and adaptive
  mesh-refinement techniques}.
\newblock Teubner Skripten zur Numerik, 1996.

\end{thebibliography}
\end{document}